\documentclass[notitlepage,twoside,a4paper]{amsart}
\usepackage{amssymb}
\usepackage[totalwidth=15cm,totalheight=23cm,marginparwidth=40pt]{geometry}

\newcommand{\II}{I\hspace{-0.1cm}I}
\newcommand{\III}{I\hspace{-0.1cm}I\hspace{-0.1cm}I}

\newcommand{\ula}{\mathbf S}

\DeclareMathOperator{\ext}{ext}

\newtheorem{theorem}{\rm\bf Theorem}[section]

\newtheorem{definition}[theorem]{\rm\bf Definition}

\newtheorem{question}[theorem]{\rm\bf Question}

\newcommand{\C}{{\mathbb C}}
\newcommand{\CP}{{\mathbb CP}}

\newcommand{\DD}{{\mathbb D}}
\newcommand{\HH}{{\mathbb H}}
\newcommand{\R}{{\mathbb R}}

\newcommand{\cCC}{{\mathcal{CC}}}

\newcommand{\cGH}{{\mathcal {GH}}}
\newcommand{\cH}{{\mathcal H}}

\newcommand{\cS}{{\mathcal S}}

\newcommand{\cT}{{\mathcal T}}
\newcommand{\cML}{{\mathcal{ML}}}

\newcommand{\PSL}{\rm{PSL}}

\newcommand{\tr}{\rm{tr}}

\title{Volume functions and boundary data of 3-dimensional hyperbolic manifolds}

\author{Jean-Marc Schlenker}
\thanks{Jean-Marc Schlenker:
University of Luxembourg, Department of Mathematics, 
Maison du nombre, 6 avenue de la Fonte,
L-4364 Esch-sur-Alzette, Luxembourg}
\email{jean-marc.schlenker@uni.lu}

\date{\today}

\begin{document}

\maketitle

\begin{abstract}
  We review recent progress on two closely related sets of questions concerning convex co-compact hyperbolic manifolds, or convex domains in those manifolds, such as their convex core. The first set of questions is to what extent the hyperbolic metric on such a manifold is uniquely determined by either of two possible geometric data on their boundary. The second aspect is the ``volume'' associated to such a manifold, such as the renormalized volume of a convex co-compact hyperbolic manifold. The relation between the two is provided by the first variation of the volume functions, which involves the two kinds of boundary data as ``conjugate'' variables.

  While progress has recently been made on some questions, others remain open. New connections have recently emerged, with physics (and in particular the AdS/CFT correspondence) as well as with probability theory (the Loewner energy). 
\end{abstract}

\tableofcontents

\section{Introduction}
\label{sc:intro}
We consider here a 3-dimensional manifold $M$ which is the interior of a compact manifold with boundary. Among the simplest examples are the product of a closed surface of genus at least $2$ by an interval, or the interior of a ``bretzel'' (a solid torus of genus $2$ or higher).

Such a manifold often admits a complete hyperbolic metric. We will be concerned here primarily with metrics which are {\em convex co-compact} in the following sense.

\begin{definition}
  Let $(M,g)$ be a complete hyperbolic manifold. A subset $C\subset M$ is {\em geodesically convex} if any geodesic segment in $M$ with endpoints in $C$ is contained in $C$. $(M,g)$ is {\em convex co-compact} if $M$ contains a non-empty, compact, geodesically convex subset. It is {\em geometrically finite} if it contains a non-empty geodesically convex subset of finite volume. 
\end{definition}

A convex co-compact hyperbolic metric on $M$ induces a complex structure on its ideal boundary $\partial_\infty M$, defining a map $\partial$ from $\cCC(M)$, the space of convex co-compact hyperbolic metrics on $M$, considered up to isotopy, to $\cT_{\partial M}$, the Teichm\"uller space of $\partial M$. The classical Bers Simultaneous Uniformization Theorem, and its extension to convex co-compact manifolds (see Section \ref{ssc:ahlfors-bers}) affirms that this map is a biholomorphism.

We will also be interested in geodesically convex subsets of a convex co-compact hyperbolic manifold $M$. An important example is the {\em convex core} $C(M)$ of $M$, defined as its smallest non-empty geodesically convex subset.

A general, but still partially conjectural and incomplete, picture emerges.
\begin{enumerate}
\item The hyperbolic structure on $M$ is uniquely determined by either of two geometric data that can be considered on its boundary. One is some kind of ``induced metric'' data (depending on the properties of $M$), while the other is some kind of ``bending data''.
\item Those two data are {\em dual}, and they are also {\em conjugate variables} (in a symplectic sense). 
\item One can also associate to $M$ a ``volume'', which can, depending on the properties of $M$, be its usual volume, or its {\em renormalized volume} $V_R$.
\item The first variation of this volume can be expressed in a simple way in terms of the boundary data. (The Legendre transform of this volume with respect to the boundary data, called the {\em dual volume} can also be useful in some situations.)
\item The volume of convex co-compact manifolds is closely related to quantities defined on the boundary, in particular the Liouville functional and the Loewner energy.
\item The gradient flow of the (renormalized) volume is a powerful tool to prove properties of $V_R$ over $\cCC(M)$. 
\end{enumerate}


Points (4) and (5) have a strong relation to the notion of {\em holography} in physics, and specifically with the AdS/CFT correspondence, as introduced in \cite{maldacena}. In a precise version \cite{witten}, the correspondence states the equality between the partition function of a certain CFT (depending on an integer parameter) on a $d$-dimensional manifold $X$, and the sum, over all $d+1$-dimensional manifolds $M$ with boundary $X$, of a function of the action of a superstring theory on $M$. Points (4) and (5) might appear as very simple consequences, for $d=2$ and in a limit case where many features of the full correspondence disappear. However even in this simplified limit, the AdS/CFT correspondence suggests new, simple conjectural statements on hyperbolic 3-manifolds (see Section \ref{ssc:adscft}).


\subsection*{Acknowledgement} The author would like to thank LP for a careful re-reading and for helpful comments.

\section{Convex domains in hyperbolic manifolds}

In this section we consider convex domains in hyperbolic manifolds, focusing on the convex cores of convex co-compact hyperbolic manifolds.

\subsection{Convex domains in $\HH^3$}

When the ambient manifold is simply $\HH^3$, compact, geodesically convex subsets are homeomorphic to balls. There are two kinds of results that describe them in terms of some induced data on their boundary. In each case, the statements were discovered first for polyhedra, then for domains with smooth boundary.

\subsubsection*{Induced metrics}

The first result describes a convex domain in terms of the induced metric on the boundary. When $P\subset \HH^3$ is a compact (convex) polyhedron, the induced metric on its boundary is a hyperbolic metric with cone singularities of cone angles at most $2\pi$. Conversely, Alexandrov \cite{alex} proved in 1948 that any hyperbolic metric on the sphere with at least three cone singularities of cone angles less than $2\pi$ can be realized uniquely as the induced metric on the boundary of a (possibly degenerate) polyhedron in $\HH^3$. 

A similar result was then proved by Pogorelov \cite{Po} for compact, convex domains with smooth, strongly convex boundary -- strongly convex here means that the principal curvatures are non-zero. The induced metric is then a smooth metric of curvature $K>-1$ on the boundary, and any such metric can be uniquely obtained. Pogorelov proved in fact that {\em no} curvature assumption is necessary: any non-smooth metric on $S^2$, with curvature larger than $-1$ in the sense of Alexandrov, can be realized on the boundary of a convex subset of $\HH^3$. 

\subsubsection*{Bending data}\label{sssc:domain-bending}

There is another possible way to describe convex subsets of $\HH^3$ from boundary data is in terms of ``bending data'' rather than induced metric. The first instance concerns compact or ideal polyhedra in $\HH^3$, for which the bending data is, to some extent, determined by the dihedral angles. Andreev first characterized the possible dihedral angles of compact \cite{andreev} and ideal \cite{andreev-ideal} polyhedra in $\HH^3$ with acute angles, and showed that those angles uniquely characterize a polyhedron.

Hodgson and Rivin \cite{HR} then showed that, for compact polyhedra, a complete picture can be given, in terms of the {\em dual metric}, rather than dihedral angles, of a polyhedron $P\subset \HH^3$. This dual metric can be constructed by gluing isometrically the spherical polygons defined at each vertex as the set of unit vectors orthogonal to the oriented support planes of $P$ -- when two vertices are connected by an edge $e$, the corresponding polygons each have an edge of length equal to the exterior dihedral angle of $P$ at $e$, and can therefore be glued along those edges. The metric space constructed in this way is a spherical metric with cone singularities with cone angles larger than $2\pi$ on $S^2$, and its closed geodesics have length larger than $2\pi$. It is called the {\em dual metric} of $P$. (Another interpretation of this dual metric can be found in Section \ref{ssc:ds}.)

\begin{theorem}[Hodgson, Rivin] \label{tm:HR}
  Let $h^*$ be a spherical metric with cone singularities on $S^2$, with cone angles larger than $2\pi$ and closed geodesics of length larger than $2\pi$. There is a unique compact polyhedron in $\HH^3$ with dual metric $h^*$.  
\end{theorem}

Rivin \cite{rivin-annals} also gave a complete characterization of ideal hyperbolic polyhedra in terms of their dihedral angles.

For domains with smooth, strongly convex boundary, the dual metric is replaced by another bending data, the third fundamental form of the boundary, defined as $\III(X,Y)=I(BX,BY)$, where $B$ denotes the shape operator of the boundary. It has curvature $K<1$ and closed geodesics of length larger than $2\pi$. Any such smooth metric on $S^2$ can be obtained uniquely as the third fundamental form of a convex domain with smooth, strictly convex boundary.

\subsubsection*{The Schl\"afli formula}

Denote by $V$ the volume of a compact polyhedron $P\subset\HH^3$. Let $\dot P$ be a first-order deformation of $P$, that is, the data of a vector for each of the vertices of $P$, such that the faces of $P$ remain planar under the deformation (the combinatorics of $P$ does not change, and no face is ``broken''). Let $E$ be the set of edges of $P$, and let $\theta_e$ be the exterior dihedral angle of $P$ at edge $e$. Schl\"afli \cite{schlafli1858,milnor-schlafli} established a simple first-order variation formula for the volume of $P$ in terms of its boundary invariant:
\begin{equation}
  \label{eq:schlafli}
  \dot V = \frac 12 \sum_e l_e \dot\theta_e~. 
\end{equation}
Here the sum is over the edges of $P$, $l_e$ is the edge length of edge $e$, while $\theta_e$ is its exterior dihedral angle.

This formula shows the interplay between the two types of boundary data, the induced metric (appearing as edge lengths) and the bending data (appearing as dihedral angles) in the first-order variation of the volume.




There is an analog of the Schl\"afli formula, valid for deformations of a convex domain $N$ in $\HH^3$, or in fact in hyperbolic manifolds, see \cite{sem-era}.
\begin{equation}
  \label{eq:sem}
  \dot{V(n)} = -\frac 12\int_{\partial N} \dot H + \frac 12\langle \dot I, \II\rangle_I da_I~, 
\end{equation}
where $I$ denotes the induced metric on $\partial N$, $\II$ denotes its second fundamental form, and $H$ is its mean curvature (defined as $H=\tr_I\II$). This equation can also be written as
\begin{equation}
  \label{eq:HE}
  \dot V + \frac 12 \dot \cH = - \frac 14 \int_{\partial N} \langle \dot I,\II_0\rangle_I da_I~,
\end{equation}
where $\cH$ is the integral mean curvature of $\partial N$, $\II_0$ is the traceless part of the second fundamental form of $\partial N$. In this form, the variational formula appears closely related to the variational formula for the Hilbert-Einstein formula with its natural boundary term (called the Gibbons-Hawking-York term in the physics literature).





\subsection{The convex core}
\label{ssc:core}

\subsubsection*{Definition}

By definition, a convex co-compact hyperbolic manifold $M$ contains a non-empty, compact, geodesically convex subset. It follows from the definition that the intersection between two geodesically convex subsets $K_1$ and $K_2$ is geodesically convex, and it can also be proved that $K_1\cap K_2\neq \emptyset$ if $K_1\neq \emptyset$ and $K_2\neq \emptyset$. Therefore, $M$ contains a smallest non-empty geodesically convex subset, its {\em convex core} $C(M)$. It can be proved (see \cite[Section 8.8]{thurston-notes}) that $M\setminus C(M)$ is the disjoint union of a finite set of ``hyperbolic ends'', each homeomorphic to the product of a surface of genus at least $2$ with $[0,\infty)$. So $\partial C(M)$ is homeomorphic to $\partial_\infty M$.

The simplest example is provided by quasifuchsian manifolds, which are the convex co-compact manifolds homeomorphic to $S\times \R$, where $S$ is a closed surface of genus at least $2$.

A good reason to study the convex cores of convex co-compact hyperbolic manifolds is that they provide -- even in simple cases, for quasifuchsian manifolds homeomorphic to $S\times \R$ -- a ``bridge'' between two points of view on Teichm\"uller theory:
\begin{itemize}
\item the ``hyperbolic'' point of view, where points in $\cT_S$ correspond to hyperbolic structures on $S$, cotangent vectors correspond to measured laminations, etc,
\item the ``complex'' point of view, where points in $\cT_S$ correspond to complex structures, cotangent vectors to holomorphic quadratic differentials, etc. 
\end{itemize}
While the data on $\partial_\infty M$ are mostly related to the ``complex'' Teichm\"uller theory, the data on $\partial C(M)$ is mostly in terms of hyperbolic geometry.

\subsubsection*{The structure of the boundary}

It follows from its definition that $C(M)$ is a minimal geodesically convex subset, so it can not have any extreme point, since otherwise those points could be ``cut off'' by a plane to make $C(M)$ smaller (but still geodesically convex). From this starting point, Thurston \cite[Section 8.8]{thurston-notes} proved that each connected component of $\partial C(M)$ is a {\em pleated surface}: its induced metric is hyperbolic (of constant curvature $-1$) and it is ``bent'' along a measured geodesic lamination. We will denote by $\cML_{\partial M}$ the space of measured geodesic laminations on $\partial M$ (see \cite{FLP,bonahon:laminations} or \cite[Chapter 8]{thurston-notes} for more on this notion). 

A convex co-compact hyperbolic metric $g\in \cCC(M)$ therefore determines two boundary data on $\partial C(M)$, similar to those mentioned above for polyhedra or convex domains with smooth boundary in $\HH^3$:
\begin{itemize}
\item a hyperbolic metric $m\in \cT_{\partial M}$, defined as the induced metric on $\partial C(M)$,
\item a measured lamination $l\in \cML_{\partial M}$, defined as the bending, or pleating, measured lamination of $\partial C(M)$.
\end{itemize}
Thurston asked whether either of these boundary data uniquely determines $g$.

\subsubsection*{The induced metric on $\partial C(M)$}

For the induced metric on $\partial C(M)$, it is known that any hyperbolic metric can be obtained -- the following result follows from work or Epstein--Marden \cite{epstein-marden}, or independently, of Labourie \cite{L4}.

\begin{theorem} \label{tm:induced}
  Let $h\in \cT_{\partial M}$ be a hyperbolic metric on $\partial M$. There exists a hyperbolic metric $g\in \cCC(M)$ such that the induced metric on the boundary of the convex core is isometric to $h$.
\end{theorem}

The uniqueness of a hyperbolic structure on $M$ such that the induced metric $m$ is prescribed, however, remains elusive.

\begin{question}[Thurston] \label{q:induced}
  Let $m\in \cT_{\partial M}$ be a hyperbolic metric. Is there a {\em unique} hyperbolic metric $g\in \cCC(M)$ such that the induced metric on the boundary of the convex core is $m$? 
\end{question}

There is an analogy between the induced metric $m$ and the conformal metric at infinity $c$ of a convex co-compact hyperbolic manifold. However, for manifolds with incompressible boundary, the two are also {\em close} in the sense that there exists a universal constant $K_0$ such that $m$ is $K_0$-quasi-conformal to $c$ -- this was first established by Sullivan \cite{sullivan:travaux}, and further refined in \cite{bishop:explicit,epstein-marden-markovic,bridgeman-canary-yarmola}. For convex co-compact manifolds with compressible boundary, the constant depends on the injectivity radius of $m$ \cite{canary:conformal}.

\subsubsection*{The measured bending lamination}

The description of a convex co-compact hyperbolic metric in terms of its measured bending lamination is better understood. Bonahon and Otal \cite{bonahon-otal} described the measured laminations that can be obtained in this manner for geometrically finite manifolds with incompressible boundary, and proved that rational laminations (those supported on a finite set of closed curves) are uniquely realized. Their result was extended by Lecuire \cite{lecuire} to manifolds with compressible boundary. The fact that a convex co-compact hyperbolic metric is uniquely determined by the measured bending lamination on the boundary of its convex core was proved in \cite{uniqueness}.

\begin{theorem}[Bonahon--Otal, Lecuire, Dular--S]\label{tm:bonahon-otal}
  Let $\bar M$ be a compact 3-manifold with non-empty boundary, with all boundary components of genus at least 2, and such that the interior $M$ of $\bar M$ admits a complete hyperbolic metric. Let $l$ be a measured lamination on $\partial \bar M$. Assume that:
  \begin{enumerate}
    \item Each closed leaf of $l$ has weight less than $\pi$,
    \item For each essential disk $D$ in $\bar M$, $i(l,\partial D)>2\pi$,
    \item There exists $\eta >0$ such that $i(\partial A,l)\geq\eta$ for each essential annulus $A$ in $\bar M$.
  \end{enumerate}
  Then there exists a unique non-Fuchsian, convex co-compact hyperbolic metric on $M$ such that $l$ is the measured pleating lamination of the boundary of the convex core of $M$.
\end{theorem}

Note that the results of Bonahon--Otal and Lecuire are in fact more general, since they apply to geometrically finite manifolds. Fuchsian manifolds are excluded from the statement since their measured bending lamination is always zero. 





\subsubsection*{The Bonahon-Schl\"afli formula}

Since it is compact, the convex core of a quasifuchsian manifold $M$ has finite volume. We denote by $V_C:\cCC(M)\to \R_{\geq 0}$ the corresponding volume function.

Bonahon \cite{bonahon-variations,bonahon} extended the Schl\"afli formula to the convex core in a rather subtle manner. He proved that under a first-order deformation of the hyperbolic structure on $M$,
\begin{equation}
  \label{eq:bonahon-schafli}
  V_C' = \frac 12 L_m(l')~. 
\end{equation}
Here $l'$ denotes the first-order variation of the measured bending lamination, which he defines as a {\em H\"older cocycle}, and $L_m(l')$ denotes its hyperbolic length, which Bonahon also defines. Thus $l'$ plays the role of the first-order variation of the dihedral angles for polyhedra.

The formula becomes simpler for the Legendre transform of the volume function, called {\em dual volume} for reasons that will appear in Section \ref{ssc:ds}. It is defined as
$$ V^*_C = V_C - \frac 12 L_m(l)~, $$
and its first-order variation is described by the dual Bonahon-Schl\"afli formula:
\begin{equation} 
  \label{eq:bs-dual}
  (V^*_C)' = -\frac 12 dL_\cdot(l)(m')~. 
\end{equation}
This variational formula is simpler since the right-hand term is simply the differential of the length of $l$ -- an analytic function over $\cT_{\partial M}$ (see \cite{kerckhoff:analytic}) applied to the tangent vector $m'$. 

\subsection{Larger geodesically convex subsets}
\label{ssc:larger}

Some of the properties described above for polyhedra in $\HH^3$, or for convex cores of convex co-compact hyperbolic manifold, have extensions to larger compact, geodesically convex domains in convex co-compact manifolds.

\begin{theorem} \label{tm:hmcb}
\begin{enumerate}
\item If $N\subset M$ is such a domain with smooth, strongly convex boundary, then the induced metric $I$ has curvature $K>-1$. Any smooth metric of curvature $K>-1$ on $\partial M$ can be obtained in a unique manner.
\item Under the same conditions, the third fundamental form of $\partial N$ has curvature $K<1$, and its closed, contractible (in $N$) geodesics have curvature $L>2\pi$. Any smooth metric on $\partial M$ satisfying those conditions can be uniquely realized in this manner.
\end{enumerate}
\end{theorem}

The existence in part (1) was proved in \cite{L4}, the uniqueness in \cite{hmcb}. In fact, any (non-smooth) metric with curvature $K>-1$ in the sense of Alexandrov can be realized on the boundary of a geodesically convex subset of $M$ equipped with a convex co-compact metric on $M$, see \cite{slutskiy:compact}. Point (2) was obtained in \cite{hmcb}.

The volume of $N$ satisfies the first-order variation formula \eqref{eq:sem} (or equivalently \eqref{eq:HE}).

It is also quite natural to consider geodesically convex subsets with polyhedral boundary in a convex co-compact hyperbolic manifold $M$ --- they are the convex hulls of finite set of points, in the sense that the are the smallest geodesically convex subsets of $M$ containing a given, finite set of points. Prosanov \cite{prosanov:dual} recently extended to this polyhedral setting the existence and uniqueness result of \cite{HR} for convex polyhedra. He also made progress on the existence and uniquess question for the induced metric \cite{prosanov:polyhedral}. 

\subsection{Domains in $\HH^3$}
\label{ssc:domains}

Let $M$ be a quasifuchsian manifold. The convex core $C(M)$ lifts in the universal cover of $M$, identified with $\HH^3$, to the convex hull in $\HH^3$ of the limit set $\Lambda_\rho$ of the holonomy representation $\rho:\pi_1M\to PSL(2,\C)$. The boundary of the convex hull of $\Lambda_\rho$, $\partial CH(\Lambda_\rho)$, is then the disjoint union of two pleated disks, each isometric to $\HH^2$, which are ``glued'' at infinity along $\Lambda_\rho$ by a map $u:\partial_\infty \HH^2\to \partial_\infty\HH^2$ which is equivariant under the holonomy representations $\rho_+, \rho_-$ of the induced metrics on the two boundary components of $\partial C(M)$, in the sense that
$$ \forall \xi \in \partial_\infty \HH^2, \forall \gamma \in \pi_1M, u(\rho_-(\gamma)\xi) = \rho_+(\gamma) u(\xi)~. $$

Theorem \ref{tm:induced} can be stated in this context as the existence, for each quasi-symmetric homeomorphism $u:\partial_\infty \HH^2\to \partial_\infty\HH^2$ which is equivariant under two Fuchsian representations of $\pi_1M$ on $\HH^2$, of a quasi-circle $\Gamma\subset \partial_\infty \HH^3$ such that the two connected components of $\partial CH(\Gamma)$, which are both isometric to $\HH^2$, are identified at infinity, along $\Gamma$, by $u$.

Similarly, part (1) of Theorem \ref{tm:hmcb} can be stated, in the special case of a quasifuchsian manifold, as follows: given two Fuchsian representations $\rho_+, \rho_-;\pi_1M\to \PSL(2, \R)$, given a quasi-symmetric homeomorphism $u:\partial\DD^2\to \partial\DD^2$ equivariant under $\rho_+, \rho_-$, and given two conformal metrics $h_+, h_-$ of curvature $K>-1$ on the disk $\DD$, there exists a unique convex subset $C\subset \HH^3$, invariant under a quasifuchsian representation $\rho:\pi_1M\to PSL(2,\C)$ and with ideal boundary the limit set $\Lambda_\rho$, such that the induced metric on the upper (resp. lower) boundary components of $\partial C$ are isometric to $h_-$ and $h_+$, respectively, with the disks identified at infinity along $\Lambda_\rho$ by $u$.

One can then ask whether the existence, and possibly the uniqueness, of $\Gamma$ still holds when the quasi-symmetric homeomorphism $u:\partial_\infty \HH^2\to \partial_\infty\HH^2$ is not supposed to be equivariant. An existence result was proved in \cite{convexhull} for metrics of constant curvature $K\in [-1,0)$, and in \cite{chen2022geometric} for metrics of variable curvature. Mesbah \cite{mesbah2024} also obtained an existence result for the corresponding problem in the anti-de Sitter space (see Section \ref{ssc:ads}).

However, at this point, only existence results are known in this ``universal'' context, and all uniqueness statements remain conjectural.

Note also that if $\Gamma\subset \partial_\infty \HH^3$ is a quasi-circle, the volume of $CH(\Gamma)$ is generally infinite, unless $\Gamma$ is assumed to be {\em Weil--Petersson} (see Section \ref{sc:loewner}).

\section{Convex co-compact hyperbolic manifolds}

We outline in this section some results and questions which are similar to those described in the previous section on the boundary of the convex core. The induced metric is replaced by the conformal structure at infinity, the volume (of dual volume) of the convex core is replaced by the renormalized volume, while the measured bending lamination of $\partial C(M)$ is replaced by a ``measured foliation at infinity'', which seems to play a similar role.

\subsection{The Ahlfors--Bers theorem}
\label{ssc:ahlfors-bers}

Given an orientable convex co-compact hyperbolic metric $g$ on $M$, the boundary at infinity $\partial_\infty M$ can be identified with the quotient of an open subset of $\partial_\infty\HH^3$ by an action of the fundamental group of $M$. This action is by elements of $PSL(2,\C)$, so that $\partial_\infty M$ is endowed with
\begin{itemize}
\item a complex projective structure $\sigma_\infty$ (see \cite{dumas-survey}),
\item a complex structure $c_\infty$, which is the underlying complex structure of $\sigma_\infty$.
\end{itemize}

The Ahlfors--Bers theorem, due to the efforts of several authors including Ahlfors, Bers, Kra, Marden, Maskit, Sullivan and Thurston (see e.g. \cite[\S 5.1, 5.2]{marden:hyperbolic} or \cite[\S 5.3]{matsuzaki-taniguchi}) asserts that the map sending $g\in \cCC(M)$ to $c_\infty\in \cT_{\partial M}$ is a biholomorphism. (Note that this holds because $g$ is considered up to isotopy.) We consider the uniqueness part of this theorem as an analog ``at infinity'' of Question \ref{q:induced}. We will see in Section \ref{ssc:foliation} what could be a analog of Theorem \ref{tm:bonahon-otal}.

\subsection{The renormalized volume of convex co-compact manifolds}

\subsubsection*{Origins}

Convex co-compact hyperbolic manifolds always have infinite volume (unless they are compact). However, in 1998, Henningson and Skenderis \cite{Skenderis} proposed the definition of a ``renormalized'' volume for conformally compact Einstein manifolds -- a notion that, in dimension 3, coincides with that of convex co-compact hyperbolic manifold. Their definition essentially reduced to disregarding diverging terms in a suitable expansion of the volume of increasing regions in the manifold. The definition was then brought into the mathematics realm, and extended to minimal submanifolds of Einstein conformally compact manifolds, by Graham and Witten \cite{graham-witten,graham}. (See also \cite[Section 2]{witten} for a clear, physically motivated explanation.)

For conformally compact Einstein manifolds of even dimension, the renormalized volume is well-defined. In odd dimension, however, it depends on the choice of a metric in the conformal class at infinity. For 3-dimensional manifolds, however, there is a canonical choice of a metric in each conformal class -- the unique hyperbolic metric, since each boundary component has genus at least $2$.

\subsubsection*{Definition by $K$-surfaces and $W$-volume}
\label{sssc:dfVR}

A simple definition of the renormalized volume is given by Mazzoli \cite{mazzoli2019-1}. It is based on a result of Labourie \cite{L6}: given a convex co-compact hyperbolic manifold $M$, each end of $M$ -- that is, each connected component of $M\setminus C(M)$ -- admits a unique foliation by surfaces of constant curvature $K$, with $K\in [-1,0)$. Let $N_K$ be the geodesically convex domain in $M$ bounded by the $K$-surfaces in each end.

We consider a modified version of the volume of $N_K$, its $W$-volume \cite{volume}, defined as
\begin{equation}
  \label{eq:W}
  W(N_K) = V(N_K) - \frac 14\int_{\partial N_K} Hda_I~,
\end{equation}
where $H=\tr(B)$ is the trace of the shape operator of $\partial N_K$. 

The renormalized volume can then be defined \cite[Theorem C]{mazzoli2019-1} as
$$ V_R(M) = \lim_{K\to 0^-} W(N_K) - \pi |\chi(\partial M)|\rm{arctanh}(\sqrt{K+1})~. $$

There is another, slightly more complicated, but perhaps more natural way to define the renormalized volume. Any compact, geodesically convex subset $N\subset M$ defines a metric $h_N$ in the conformal class of $\partial_\infty M$, obtained by suitably scaling the induced metric on the surface at constant distance $r$ from $N$ as $r\to \infty$. One would like to say that there is a unique choice of $N$ such that $h_N$ is hyperbolic. This is equivalent to taking for $N$ the domain bounded by the Epstein surface of the hyperbolic metric at infinity (defined in \cite{epstein:envelopes}). (In practice, this Epstein surface might in some cases not be embedded or even immersed, but this difficulty can be avoided by suitably scaling the hyperbolic metric at infinity, see \cite{volume}.)

\subsubsection*{Relations to the volume of the convex core}

Some of the applications of the renormalized volume mentioned below use the fact that the renormalized volume of a convex co-compact hyperbolic manifold $M$ is coarsely equivalent, and even within bounded additive constants, from the volume of the convex core.

For manifolds with incompressible boundary, this was proved in \cite{compare} for quasifuchsian manifolds, and extended to general convex co-compact manifolds by Bridgeman and Canary \cite{bridgeman-canary:renormalized}. Specifically, they prove that
\begin{equation}
  \label{eq:bridgeman-canary}
  V_C(M)-10|\chi(\partial M)|\leq V_R(M)\leq V_C(M)~. 
\end{equation}

For manifolds with compressible boundary, Bridgeman and Canary \cite{bridgeman-canary:renormalized} show that $V_R(M)-V_C(M)$ is bounded from below and from above by terms which are linear in the logarithm of the length of the shortest geodesic on $\partial  C(M)$ which bounds a disk in in $C(M)$.

However, even for manifolds with compressible boundary, a close relation exists between the renormalized volume of $M$ and the $W$-volume of its convex core. The $W$-volume of the convex core is equal (as can be seen from \eqref{eq:W} by a limiting argument) to
$$ W(C(M))=V_C(M) - \frac 14 L_m(l)~. $$
It is then proved in \cite[Lemma A.6]{averages} that
$$ W(C(M)) -C(\partial M) \leq V_R(M) \leq W(C(M)) + \frac{\pi\log(2)}2 |\chi(\partial M)|~, $$
where $C(\partial M)$ only depends on the topology of $\partial M$. 

\subsection{The variational formula of $V_R$}

A key property of the renormalized volume is its first variation property. We provide it here in Section \ref{sssc:first}, after some definitions. We will see in Section \ref{ssc:foliation} that it is in fact closely analogous to the (dual) Bonahon-Schl\"afli formula \eqref{eq:bs-dual}.

\subsubsection*{The Schwarzian derivative at infinity}

Let $M$ be a quasifuchsian hyperbolic manifold. As seen in Section \ref{ssc:ahlfors-bers}, its asymptotic boundary is equipped with a complex projective structure $\sigma_\infty$, to which can be associated a holomorphic quadratic differential $q$ on $(\partial_\infty M, c_\infty)$. Specifically, each connected component of $\partial_\infty M$ can be written as the quotient of a domain $\Omega\subset \partial_\infty \HH^3\simeq \C P^1$ by a surface group acting by elements of $PSL(2,\C)$. Then $q$ is equal to the Schwarzian derivative of the (inverse) uniformization map from $\Omega$ to the disk. We call $q$ the {\em Schwarzian derivative at infinity} or {\em Schwarzian at infinity} of $M$.

A special case is when $\partial C(M)$ is totally geodesic. Then each connected component of $\partial_\infty M$ is the quotient of a disk in $\CP^1$ by an action by M\"obius transformations, and therefore $q=0$. In fact, $q=0$ exactly when $\partial C(M)$ is totally geodesic. 

This holomorphic quadratic differential can be used to bound the geometric invariants of $\partial C(M)$. Specifically, \cite{bridgeman-brock-bromberg} contains among other results
\begin{itemize}
\item a bound on the Lipschitz constant of the retraction map $\partial_\infty M$ to $\partial C(M)$ in terms of the $L^\infty$ norm of $q$,
\item a bound on the hyperbolic length of the measured bending lamination $l$ of $\partial C(M)$ by
  $$ L_m(l)\leq 4|\chi(\partial M)| \|q\|_\infty~, $$
\item a bound on $V_R(M)-V_C(M)$ in terms of $\|q\|_2$ when the convex co-compact structure is close to that for which the convex core has geodesic boundary, specifically,
$$ V_C(M) - |\chi(\partial M)|G(\| q\|_2) \leq V_R(M)\leq V_C(M)~, $$
with $G(t)\sim t^{1/5}$ at $t=0$.
\end{itemize}

\subsubsection*{A first variation formula}
\label{sssc:first}

A first-order variation of the hyperbolic structure on $M$ is associated, through the Ahlfors--Bers Theorem, to a first-order variation $\dot c_\infty$ of the complex structure at infinity. This variation $\dot c_\infty$ corresponds to a Beltrami differential on $(\partial_\infty M, c_\infty)$, which can be paired with the holomorphic quadratic differential $q$. The first-order variation of the renormalized volume is then
\begin{equation}
  \label{eq:dVR}
  \dot V_R = \mathrm{Re}\left( \langle \dot c_\infty, q\rangle\right)~,
\end{equation}
where $\langle, \rangle$ denotes the duality pairing.

A clear analogy with the variational properties of the dual volume of the convex core appears (see \eqref{eq:bs-dual}): in both cases, the right-hand term of the formula is equal to the restriction, to the variation of the boundary data ($(m,l$) in the first case, $(c_\infty,q)$ in the second), of the Liouville form of an ambient symplectic form. This analogy is even stronger after identifying $\cT_{\partial M}\times \cML_{\partial M}$ with $T^*\cT_{\partial M}$, through the map
$$ (m,l)\mapsto (m, dL_m(l))~, $$
where $dL_m(l)$ denotes the differential at $m$ of the (analytic) function on $\cT_{\partial M}$ sending $m$ to the hyperbolic length $L_m(l)$ of $l$ for $m$.

There is, however, an even stronger analogy, which will be described in Section \ref{ssc:foliation}.

\subsubsection*{Relations to the Liouville functional}

In many situations, such as quasifuchsian or Schottky manifolds, the renormalized volume corresponds to the Liouville functional defined previously, in particular by Takhtajan--Zograf \cite{TZ-schottky} and Takhtajan--Teo \cite{takhtajan-teo}.

The correspondence between the Liouville functional and the renormalized volume was discovered by Krasnov \cite{Holography}. It can be considered as another instance of holography, with a 2-dimensional quantity defined on the conformal boundary at infinity --- the Liouville functional --- equal to a 3-dimensional quantity in the ``bulk'' --- the renormalized volume.

\subsubsection*{Gradient bounds}

A consequence of the variational formula \eqref{eq:dVR} is that the Schwarzian at infinity $q$ is the differential of $V_R$, considered (through the Ahlfors--Bers Theorem) as a function defined on $\cT_{\partial M}$.

When $M$ has incompressible boundary, $q$ is the Schwarzian of a univalent map defined on a simply connected domain in $\CP^1$, and it is thus, by the classical Kraus-Nehari estimate \cite{nehari-bams}, uniformly bounded relative to the hyperbolic metric on $\partial_\infty M$. Specifically, this argument is used in \cite{compare,bridgeman-brock-bromberg:gradient} to prove that
\begin{equation}
  \label{eq:gradient}
  \| dV_R\|_{WP} \leq 3\sqrt{\frac{\pi |\chi(\partial M)|}2}~. 
\end{equation}
We will see in Section \ref{sssc:brock} that this upper bound has some applications.

When $M$ has compressible boundary, the upper bound on the gradient depends on the length of the shortest geodesic in $\partial_\infty M$ bounding a disk in $M$. Section \ref{sssc:adapted} presents an adapted version of the renormalized volume of manifolds with compressible boundary which has a Weil--Petersson gradient that is uniformly bounded.

\subsection{Positivity of the renormalized volume and the gradient flow}

The term ``volume'' suggests a quantity which is always positive. For the renormalized volume, this is far from obvious, although it is indeed the case for manifolds with incompressible boundary (and with a suitable normalization of the renormalized volume). Specifically, Bridgeman, Brock and Bromberg \cite{bridgeman-brock-bromberg:gradient} proved the following result. They consider a relatively acylindrical 3-dimensional manifold $(M,S)$, and a deformation space $\cCC(M, S; X)$.  

\begin{theorem} \label{tm:BBB}
  The function $V_R: \cCC(M, S; X)\to \R$ admits a unique minimum, attained at the unique point where $C(M)$ has totally geodesic boundary. 
\end{theorem}

To prove this result, the authors use the (inverse, i.e. decreasing) gradient flow of the renormalized volume, introduced in \cite{bridgeman-brock-bromberg}, and they show that a suitably adapted version of the flow (with surgeries close to the Weil--Petersson completion of $\cT_{\partial M}$) always converge to the unique critical point of $V_R$ on the deformation space $\cCC(M, S; X)$, namely, the hyperbolic structure for which the boundary of the convex core is totally geodesic.

This approach was then refined in \cite{bridgeman-bromberg-pallete,bridgeman-bromberg-pallete:convergence} where the authors show that the flow converges to this structure with totally geodesic boundary, without the need of any surgery.

For manifolds with compressible boundary, however, the renormalized volume is {\em not} positive -- in fact it goes to $-\infty$ upon pinching a closed, compressible curve (see \cite{bridgeman-canary:renormalized,pallete:schottky,averages}). This divergence along the Weil--Petersson metric boundary of $\cT_{\partial M}$ also implies that the Weil--Petersson gradient of $V_R$ is unbounded. An attempt at correcting those defects is described in Section \ref{sssc:adapted}.

\subsection{The measured foliation at infinity and the Schl\"afli formula}
\label{ssc:foliation}

The analogy between the variation formulas \eqref{eq:bs-dual} and \eqref{eq:dVR} can be made more striking through the following definition.

\begin{definition} \label{df:foliation}
  Let $M$ be a convex co-compact hyperbolic manifold. Its {\em measured foliation at infinity} is the horizontal measured foliation of the Schwarzian differential at infinity $q$.
\end{definition}

This definition rests on the existence of a pair of measured foliations associated to any holomorphic quadratic differential on a Riemann surface (see e.g. \cite{FLP}).  There is a close analogy between measured laminations on a hyperbolic surface and measured foliations on a Riemann surface (see also \cite{FLP}). In this analogy, the hyperbolic length $L_m(l)$ of a measured lamination with respect to a given hyperbolic metric $m$ is replaced by the extremal length of a measured foliation $f$ for a given complex structure $c$. We will denote by $\cML_S$ the space of measured foliations on a closed surface $S$.

The following variational formula for the renormalized volume of a convex co-compact hyperbolic manifold is a relatively direct consequence of \eqref{eq:dVR}, see \cite[Theorem 1.7]{volumes}.
\begin{equation}
  \label{eq:dVR-foliations}
  V_R' = -\frac 12 d\ext(f)(\dot c_\infty)~,
\end{equation}
where $f$ is the foliation at infinity of $M$, $\ext(f)$ is its extremal length, considered as a function over $\cT_{\partial M}$, and $\dot c_\infty$ is the first-order variation of the complex structure at infinity, seen as a tangent vector to $\cT_{\partial M}$. 

\subsubsection*{Prescribing the measured foliation at infinity}

There is obviously a close analogy between the dual Bonahon-Schl\"afli formula \eqref{eq:bs-dual} and the variational formula \eqref{eq:dVR-foliations}. It suggests to expand this analogy between the measured foliation at infinity $f$ and the measured bending lamination on the boundary of the convex core $l$. Looking for an analog of Theorem \ref{tm:bonahon-otal} leads to a series of questions.

\begin{question}
  Let $\bar M$ be a compact 3-manifold with non-empty boundary, with all boundary components of genus at least 2, and such that the interior $M$ of $\bar M$ admits a complete hyperbolic metric.
  \begin{enumerate}
  \item Can one describe the space of measured foliations on $\partial M$ which can be obtained as the measured foliations at infinity of a convex co-compact (resp. geometrically finite) manifold?
  \item For instance, if $M=S\times \R$, where $S$ is a closed surface of genus at least $2$ equipped with a quasifuchsian metric, do the measured foliations on the two connected component of the boundary at infinity of $M$ fill $S$? Is any other condition necessary?
  \item Does its measured foliation at infinity uniquely determine a convex co-compact (or geometrically finite) hyperbolic metric on $M$? 
  \end{enumerate}
\end{question}

Recently, Choudhury and Markovi\'c \cite{choudhury-markovic} have found a partial answer to the second question. They proved that if $f_-(M), f_+(M)\in \cML_S$ denote the measured foliations at infinity on the two connected components of the boundary at infinity of a quasifuchsian hyperbolic manifold $M$, then $(f_-(M), f_+(M))$ is filling if $M$ is close enough to being Fuchsian. They also proved that if $(F_-, F_+)$ is any filling pair of measured foliations on $S$, and if $t>0$ is small enough, then the pair $(tF_-, tF_+)$ can be realized as the measured foliations at infinity of a quasifuchsian structure on $S\times \R$, which is unique near the Fuchsian locus. 

\subsection{Applications of the variational formula}

\subsubsection*{The renormalized volume as K\"ahler potential}

A first consequence of the variational formula for the renormalized volume is that $V_R$ is a K\"ahler potential for the Weil--Petersson metric on $\cT_{\partial M}$. This essentially follows from the work of Takhtajan--Zograf \cite{TZ-schottky} and Takhtajan--Teo \cite{takhtajan-teo}, see also \cite{guillarmou-moroianu-rochon} for geometrically finite manifolds.

\subsubsection*{McMullen's quasifuchsian and Kleinian reciprocity}

McMullen's quasifuchsian reciprocity \cite{McMullen} is a relation between the variations of the data at infinity --- conformal metric and Schwarzian derivative at infinity --- under a first-order deformation of a quasifuchsian hyperbolic manifold.

Let $S$ be a closed surface of genus at least $2$, and let $M=S\times \R$. Let $\cCC(M)$ be the space of quasifuchsian hyperbolic metrics on $M$, with ideal boundary components denoted by $\partial_-M$ and $\partial_+M$. As seen in Section \ref{ssc:ahlfors-bers}, each quasifuchsian metric $g\in \cCC(M)$ induces a complex projective structure $\sigma_\pm$ on $\partial_\pm M$, and therefore a pair $(c_\pm, q_\pm)\in T^*\cT_{\partial_\pm M}$, on each boundary component.

A quasifuchsian metric $g\in \cCC(M)$ is uniquely determined by $c_-, c_+$ through the Bers Simultaneous Uniformization theorem, so a first-order variation $\dot c_-$ of $c_-$ determines a first-order variation of $g$, and so does a first-order variation of $\dot c_+$ of $c_+$. They in turn determine first-order variations $\dot q_+$ of $q_+$ and $\dot q_-$ of $q_-$. This construction defines two maps
\begin{equation}
  \label{eq:phi+}
  \begin{array}[h]{cccc}
    \phi_+ & T\cT_{\partial_+M} & \to & T^*\cT_{\partial_-M} \\
    & \dot c_+ & \mapsto & \dot q_-~, 
  \end{array}
\end{equation}
\begin{equation}
  \label{eq:phi-}
  \begin{array}[h]{cccc}
    \phi_- & T\cT_{\partial_-M} & \to & T^*\cT_{\partial_+M} \\
    & \dot c_- & \mapsto & \dot q_+~. 
  \end{array}
\end{equation}

\begin{theorem}[McMullen] \label{tm:quasifuchsian}
  The maps $\phi_+$ and $\phi_-$ are adjoint.
\end{theorem}

This property extends to convex co-compact manifolds (see \cite[Appendix]{McMullen}). Given such a convex co-compact manifold $M$, $\partial_\infty M$ is equipped with a pair $(c,q)\in T^*\cT_{\partial M}$, this defines a map $\phi:\cCC(M)\to T^*\cT_{\partial M}$. Moreover $T^*\cT_{\partial M}$ is equipped with its cotangent symplectic form $\omega$.

\begin{theorem}[McMullen] \label{tm:kleinian}
  $\phi(\cCC(M))$ is Lagrangian in $T^*\cT_{\partial M}$.
\end{theorem}

Theorem \ref{tm:quasifuchsian} follows directly from this statement, using the definition of the cotangent symplectic form on $T^*\cT_{\partial M}$. The proof of Theorem \ref{tm:kleinian} follows directly from the variational formula \eqref{eq:dVR}. Indeed, \eqref{eq:dVR} can be interpreted as the fact that, on $\cCC(M)$,
$$ dV_R = \phi^*\lambda~, $$
where $\lambda$ is the Liouville form of $\omega$ on $T^*\cT_{\partial M}$. As a consequence,
$$ \phi^*\omega = \phi^*(d\lambda)=d\phi^*\lambda = d(dV_R)=0~, $$
so $\phi(\cCC(M))$ is Lagrangian in $(T^*\cT_{\partial M}, \omega)$.

\subsubsection*{Volume vs Weil--Petersson distance}
\label{sssc:brock}

Brock \cite{brock:2003} discovered a remarkable connection between the volume of the convex core of quasifuchsian manifolds and the Weil--Petersson distance between the conformal structures on the two connected components of its boundary at infinity. If $M\in \cCC(S\times \R)$ and $c_-, c_+$ are the conformal structures at infinity, then
$$  \frac 1A_g d_{WP}(c_-, c_+)-B_g\leq V_C(M) \leq A_gd_{WP}(c_-, c_+)+B_g~, $$
for constants $A_g,B_g>0$ depending only on the genus of $S$.

The upper bound on the Weil--Petersson gradient of $V_R$, together with the upper bound on the difference between $V_C$ and $V_R$, leads to an explicit upper bound. Since
$$ V_R(M) \leq 3\sqrt{\pi (g-1)} d_{WP}(c_-,c_+)~, $$
it follows from \eqref{eq:bridgeman-canary} that
\begin{equation}
  \label{eq:explicit}
  V_C(M)\leq 3\sqrt{\pi (g-1)} d_{WP}(c_-,c_+) + 20|\chi(S)|~. 
\end{equation}

Let $\phi:S\to S$ be a diffeomorophism. The mapping torus$M_\phi$ is the 3-dimensional manifold obtained from $S\times [0,1]$ by gluing $S\times \{ 0\}$ to $S\times \{ 1\}$ according to $\phi$. Thurston \cite[Chapter 7]{thurston-notes}, \cite{otal-hyperbolisation} proved that if $\phi$ is pseudo-Anosov, then $M_\phi$ is hyperbolic. The explicit upper bound \eqref{eq:explicit} can be used to bound the hyperbolic volume of $M_\phi$ in terms of the Weil--Petersson translation length, or entropy, of $\phi$. Specifically, it is proved in \cite{kojima-mcshane,brock-bromberg:inflexibility2} that
$$ V(M_\phi) \leq 3\sqrt{\pi(g-1)} \| \phi\|_{WP}~, $$
where $\| \phi\|_{WP}$ is the Weil--Petersson translation length of $\phi$.



\subsection{Convex co-compact manifolds with compressible boundary}

\subsubsection*{An adapted renormalized volume}
\label{sssc:adapted}

Consider a 3-dimensional manifold $M$ with incompressible boundary, which admits a convex co-compact hyperbolic metric. Since $\cCC(M)$ is biholomorphic to $\cT_{\partial M}$, the renormalized volume can be considered as function $V_R: \cT_{\partial M}\to \R$. As such, it has a number of pleasant properties, already mentioned above.
\begin{enumerate}
\item When $M$ is acylindrical, $V_R$ is bounded from below (and positive, if the right normalization is chosen in the definition of $V_R$),
\item $V_R$ is within a bounded additive constant (depending on the topology of the boundary) from $V_C$, the volume of the convex core,
\item the Weil--Petersson gradient of $V_R$ is uniformly bounded in the $L^\infty$ and the Weil--Petersson norm on $\cT_{\partial M}$,
\end{enumerate}

However those three properties fail for manifolds with compressible boundary. The main reason for this failure is that, when one pinches a closed, compressible curve $\gamma$ in $\partial M$, $V_R\to -\infty$.

To get over those limitations, Giovannini \cite{giovannini:adapted} introduced an {\em adapted} renormalized volume for convex co-compact manifolds with compressible boundary, as follows.

\begin{definition} \label{df:adapted}
  Let $M$ be a 3-dimensional manifold which admits a convex co-compact hyperbolic manifold. Its {\em adapted renormalized volume} is defined as the function
  $$ \tilde{V_R}: \cT_{\partial M} \to \R~, $$
  defined as
  \begin{equation}
    \label{eq:VRt}
    \tilde{V_R}(X) = V_R(X) + \max_{\mu \in \Gamma^{comp}} L(X,\mu)~,
  \end{equation}
  where $\Gamma^{comp}$ is the set of multicurves composed of compressible curves, and
  $$ L(X,\mu)= \pi^3\sum_{\gamma\in \mu} \frac 1{L_X(\gamma)}~. $$
\end{definition}

She then proves that this adapted renormalized volume has, for manifolds with compressible boundary, properties which are comparable to those of the ``usual'' renormalized volume for manifolds with incompressible boundary.
\begin{enumerate}
\item It is bounded from below,
\item it is within a bounded additive constant (depending only on the topology of the boundary) from the volume of the convex core,
\item its Weil--Petersson gradient is bounded in $L^1$ and in $L^\infty$ norm on $\cT_{\partial M}$, and therefore also in Weil--Petersson norm.
\end{enumerate}

In addition, Giovannini \cite[Theorem 4.29]{giovannini:adapted} shows that the adapted renormalized volume behaves pleasantly under pinching of a closed contractible curve (or multicurves). Suppose for instance that $(X_t)_{t\in (0,1)}$ is a one-parameter family of complex structures on $\partial M$ converging in the Weil--Petersson completion of $\cT_{\partial M}$, as $t\to 0$, to a limit $x_0$ in which a simple compressible closed curve $\gamma$ is pinched. Then the convex co-compact structure $g_t$ on $M$ corresponding to $X_t$ through the Alhfors-Bers Theorem converges in the Gromov-Hausdorff topology to a limit $g_0$, which is a convex co-compact structure on a manifold $M_0$  (which is disconnected if $\gamma$ is separating). Giovannini then proves that $\tilde{V_R}(X_t)\to V_R(M_0, g_0, \xi_0, \xi_1)$, where $V_R(M, g_0, \xi_0, \xi_1)$ denotes a renormalized volume of $(M, g_0)$ with two marked points at infinity $\xi_0, \xi_1$, defined as the ``usual'' renormalized volume of $(M, g_0)$ (see Section \ref{sssc:dfVR}) but replacing the hyperbolic metric on $\partial_\infty(M_0, g_0)$ by the hyperbolic metric of finite area on $\partial_\infty(M_0, g_0)\setminus \{ \xi_0,\xi_1\}$. 

The weakness of the adapted renormalized volume, however, is that it is only $C^{1,1}$ smooth on $\cT_{\partial M}$, since the second term on the right-hand side of \eqref{eq:VRt} is not $C^1$ at the points where several multicurves realize the maximum.

\subsubsection*{Questions stemming from the AdS/CFT correspondence}
\label{ssc:adscft}

The AdS/CFT correspondence suggests question of independent interest in hyperbolic geometry. The correspondence, in the form proposed in \cite[Section 3]{witten}, takes the form
\begin{equation}
  \label{eq:adscft}
Z_{CFT}(X) = \sum_{i} \exp(-F(M_i))~,   
\end{equation}
where $X$ is a $d$-dimensional manifold equipped with a conformal structure, the sum on the right is over all conformally compact Einstein manifolds $M_i$ with conformal boundary $X$, $Z_{CFT}(X)$ is the partition function of a CFT on $X$, and $F(M_i)$ is the action of a superstring theory on $X_i$, or possibly its product by a compact manifold. (The theories on both sides depend on an integer number $N$. Also note that the notations here differ from those used in \cite[Section 3]{witten}, where $M$ denotes the boundary conformal manifold and $X_i$ the Einstein manifolds ``filling'' it.) 

We are interested in the case where $d=2$, where the $X_i$ are the hyperbolic convex co-compact with a given (possibly disconnected) Riemann surface $M$ at infinity. In a limit where all interesting features of the superstring theory on the $M_i$ are forgotten, $F(M_i)$ is replaced simply by the renormalized volume $V_R(M_i)$. 

The dominant term on the right then corresponds to the hyperbolic manifold $M_i$, with conformal boundary $X$, of smallest renormalized volume. In case $X$ is non-connected, say the disjoint union of two connected Riemann surfaces $X^0$ and $X^1$, the left-hand side of \eqref{eq:adscft} should be the product $Z_{CFT}(X^0)Z_{CFT}(X^1)$, suggesting that the main contribution on the right-hand term should come from a disconnected filling $M=M^0\cup M^1$, with $\partial_\infty M^0=X^0$ and $\partial_\infty M^1=X^1$.

If $X^0=\overline{X^1}$, that is, $X^0$ is the same Riemann surface as $X^1$ with the orientation reversed, there exists a Fuchsian manifold filling $X^0\cup X^1$, and the heuristics above suggests that there should also be another convex co-compact hyperbolic manifold with ideal boundary $X^0$, and negative renormalized volume. It is tempting to assume that this lowest volume filling should be a Schottky manifold -- defined here as a convex co-compact hyperbolic manifold which is homeomorphic to a handlebody. This leads to the following question, attributed to Maldacena in \cite{pallete:schottky}.

\begin{question} \label{q:maldacena}
  Let $X$ be a closed, connected Riemann surface of genus at least $2$. Is there a Schottky filling of $X$ of negative renormalized volume?
\end{question}

This question remains open, although it was proved in \cite{pallete:schottky} that if $X$ is a Riemann surface of genus $g$, containing $g-1$ curves which have sufficiently short extremal length, then it does bound a Schottky manifold of negative renormalized volume. A similar result is proved in \cite{filling} under a slightly more flexible condition on the hyperbolic lengths of curves.

\subsection{The renormalized volume of Weil--Petersson quasicircles}

As when considering the determination of a domain from its boundary data in Section \ref{ssc:domains}, one can wonder whether the renormalized volume can be defined in the ``universal'' case, that is, for a quasicircle in $\partial_\infty\HH^3$ rather than for a quasifuchsian manifold. For a general quasicircle, there is so far no way to define the renormalized volume. However, recently, \cite{BBVPW} showed that a meaningful extension exists for Weil--Petersson circles.

Thanks to \cite{shen:weil-petersson}, we can use the following definition of a Weil--Petersson homeomorphism. A number of other, equivalent definition can be found in \cite{bishop:weil-petersson}.

\begin{definition}
  A quasi-symmetric homeomorphism $\phi:S^1\to S^1$ is {\em Weil--Petersson} if it is absolutely continuous with respect to the arclength measure, and $\log(\phi')\in H^{1/2}(S^1)$. A Jordan curve $\gamma\subset \partial_\infty\HH^3$ is {\em Weil--Petersson} if it is obtained from a Weil--Petersson quasi-symmetric homeomorphism by conformal welding. 
\end{definition}

We denote here by $T_0(1)$ the Universal Weil--Petersson Teichm\"uller space, defined as the space of Weil--Petersson homeomorphisms of $S^1$.

Another equivalent definition, as proved by Bishop in \cite[Theorem 1.6]{bishop:weil-petersson}, is that a quasicircle is Weil--Petersson if and only if it bounds a minimal surface in $\HH^3$ with finite integral of the squared norm of the second fundamental form (or equivalently, finite integral determinant of the shape operator). 

Recently, Bridgeman, Bromberg, Vargas Pallete and Wang \cite{BBVPW} extended to Weil--Petersson quasicircles many properties of the renormalized volume for convex co-compact hyperbolic manifolds. 

\begin{theorem} \label{tm:BBVPW}
  Let $\gamma$ be a Weil--Petersson quasicircle. Then:
  \begin{itemize}
  \item The Poincar\'e metric on the connected components of $\partial_\infty\HH^3\setminus\gamma$ still define two (possibly non-immersed) Epstein surfaces $\Sigma_\pm$ (see Section \ref{sssc:dfVR}).
  \item The signed volume $V(\gamma)$ between two Epstein surfaces is well-defined and finite.
  \item The integral mean curvatures on $\Sigma_\pm$ are also finite.
  \end{itemize}
  This makes it possible to define the renormalized volume of $\gamma$, as
  $$ V_R(\gamma) = V(\gamma) - \frac14 \int_{\Sigma_-\cup\Sigma_+} Hda~, $$
  \begin{itemize}
  \item It satisfies the variational formula
    $$ dV_R = \mathrm{Re}(\langle q, \cdot\rangle)~. $$
  \item The flow of the gradient of $V_R$ is well defined on $T_0(1)$, and all flow lines end at $0$.
  \end{itemize}
\end{theorem}

The authors also use the flow to prove that there exist universal constants $c,K>0$ such that for any $\mu\in T_0(1)$, $c(d_{WP}(\mu,0) -K)\leq V_R(\mu)$ (see \cite[Theorems 1.8 and 6.3]{BBVPW}).

\section{The renormalized volume, the Universal Liouville action, and the Loewner energy}
\label{sc:loewner}

Recently, a new connection has appeared between the renormalized volume of hyperbolic manifolds and the Loewner energy. We outline this relation in the next two section, and a (related) relation between the Schwarzian action of circle diffeomorphisms and 2-dimensional hyperbolic geometry in the following section.

Those connections between 2-dimensional quantities on Riemann surfaces (in this case, Weil--Petersson curves in $\CP^1$) and 3-dimensional quantities can again be considered as glimpses of holography, as mentioned in Section \ref{sc:intro}.

\subsection{The universal Liouville action} \label{ssc:liouville}

The universal Liouville action has been defined by Takhtajan and Teo \cite{takhtajan-teo:universal} as a K\"ahler potential for the Weil--Petersson metric on the universal Weil--Petersson Teichm\"uller space $T_0(1)$. Given a Weil--Petersson quasicircle $\gamma\subset \C$, denote by $D$ and $D^*$ the bounded and unbounded components of $\C\setminus\gamma$, and by $f:\DD\to D$ and $g:\DD^*\to D^*$ two conformal maps, with $g(\infty)=\infty$ (where $\DD\setminus \C$ is the unit disk, and $\DD^*=\C\setminus\overline{\DD}$. Then
$$ \ula(\gamma)= \int_{\DD} |f''/f'|^2 dz^2 + \int_{\DD^*}|g''/g'|^2 dz^2 + 4\pi|f'(0)/g'(\infty)|~. $$

It is proved in \cite[Theorem 1.6]{BBVPW} that if $\gamma$ is a Weil--Petersson quasicircle which in addition is $C^{5,\alpha}$ with $\alpha>0$, then $\ula(\gamma)=4V_R(\gamma)$.  

\subsection{The Loewner energy}

The Loewner energy was defined in \cite{rohde-wang} as a $L^2$ energy of the driving function of a quasicircle in $\CP^1$. It is finite for Weil--Petersson quasicircles, and it was proved in \cite{wang:equivalent} that it is equal to the universal Liouville action. It then follows from the other results of \cite{BBVPW} that it is also equal, up to a multiplicative factor, to the renormalized volume of Weil--Petersson quasicircles.

\subsection{The Schwarzian action of circle diffeomorphisms and plane hyperbolic geometry}

Ideas similar to those described in Section \ref{ssc:liouville} also appear in \cite{VPWW} in lower dimension, that is, to express quantities defined on $S^1$ in terms of quantities defined in 2-dimensional hyperbolic geometry. The main quantity considered ``at infinity'' in this 2-dimensional setting is the Schwarzian action. If $\phi:S^1\to S^1$ is a $C^3$ diffeomorphism, its Schwarzian action is defined as
$$ I(\phi) = \int_0^{2\pi} e^{2i\theta}\cS(\phi)(e^{i\theta})d\theta~, $$
where $\cS(\phi)$ is the Schwarzian derivative of $\phi$.

In \cite{VPWW}, this Schwarzian action is expressed in several ways in terms of 2-dimensional hyperbolic quantities, in particular as the (signed) hyperbolic area enclosed by the Epstein curve of the metric $\phi_*(d\theta)^2$ on $S^1$. It is also related to the Loewner energy. We refer the reader to \cite{VPWW} for more details.

\section{Lorentzian analogs}

There are close connections between questions and result described above, concerning convex co-compact hyperbolic manifolds, and questions and results concerning globally hyperbolic maximal compact (GHMC) 3-dimensional spacetimes of constant curvature.

Below we first give the definition of GHMC spacetimes of constant curvature, and the proceed to describe their relations to hyperbolic manifolds.
\begin{itemize}
\item de Sitter 3-dimensional GHMC spacetimes appear naturally as {\em projective extensions} of hyperbolic ends, and also harbour a duality with those hyperbolic ends,
\item anti-de Sitter GHMC spacetimes are {\em analoguous} to quasifuchsian hyperbolic manifolds, in the sense for instance that they also have a convex core and that many of the questions that can be asked for quasifuchsian manifolds also make sense for those GHMC AdS spacetimes.
\item Minkowski GHMC spacetimes can be interpreted as first-order deformations of Fuchsian hyperbolic manifolds (or of Fuchsian de Sitter GHMC spacetimes).  
\end{itemize}

\subsection{Globally hyperbolic maximal compact spacetimes}
\label{ssc:ghmc}

Let $M$ be a Lorentzian manifold of constant curvature. Up to scaling, we can suppose that it is modelled on either the anti-de Sitter space (of constant curvature $-1$), the Minkowski space (of constant curvature $0$), or the de Sitter space (of constant curvature $1$). We refer to \cite{mess} for the key definitions of those spaces. Recall that the de Sitter space can be seen as a quadric in the 4-dimensional Minkowski space:
$$ dS^3 = \{ x\in \R^{3,1}~|~ \langle x,x\rangle_{3,1}=1 \}~, $$
with the induced metric. The anti-de Sitter space is a similar quadric but in the flat space $R^{2,2}$ of signature $(2,2)$,
$$ AdS^3 = \{ x\in \R^{2,2}~|~ \langle x,x\rangle_{2,2}=-1 \}~, $$
while we recall for comparison that
$$ \HH^3 = \{ x\in \R^{3,1}~|~ \langle x,x\rangle_{3,1}=-1 \wedge x_0>0 \}~. $$

We say that $M$ is {\em globally hyperbolic} if it contains a Cauchy surface, that is, a space-like surface that intersects any inextendible time-like line exactly once. It is {\em globally hyperbolic compact} if this Cauchy surface is closed (compact without boundary). It is GHMC if it is globally hyperbolic compact, and maximal (in the sense of inclusion) under this condition.

\subsection{De Sitter spacetimes extending hyperbolic ends}
\label{ssc:ds}

It is useful, to understand GHMC de Sitter spacetimes, to first recall the projective duality between $\HH^3$ and the de Sitter space. Given $x\in \HH^3$, the oriented hyperplane $x^\perp$ orthogonal to $x$ is space-like in $\R^{3,1}$, and intersects $dS^3$ in a totally geodesic space-like plane. Conversely, any point in $dS^3$ is orthogonal to an oriented time-like hyperplane in $\R^{3,1}$, which intersects $\HH^3$ in an oriented totally geodesic plane. This duality has a number of pleasant properties. It also has a purely projective definition, that can be seen in the Klein (projective) model of $\HH^3$ in the unit ball $B_1$ in $\R^3$, which extends to a projective model of a hemisphere of $dS^3$ outside the ball.

If $P\subset \HH^3$, the subset of $dS^3$ of points dual to an oriented totally geodesic plane in $\HH^3$ which has $P$ on its negative side is a polyhedron $P^*$ in $dS^3$. (In the projective model of $\HH^3$ and a hemisphere of $dS^3$, it appears as the intersection with $dS^3$ of a polyhedron in $\R^3$.) The vertices of $P^*$ correspond to the 2-faces of $P$, the edges of $P^*$ to the edges of $P$, and the 2-faces of $P^*$ to the vertices of $P$. A key property used for instance in \cite{HR} is that the edge lengths of $P^*$ are the exterior dihedral angles of $P$, and the induced metric on $\partial P^*$ is the dual metric on $\partial P$. 

Let $E$ be an end of a convex co-compact hyperbolic manifold $M$, that is, a connected component of $M\setminus C(M)$. Then $E$ is bounded, on the side where it is in contact with $\partial M$, by a locally convex pleated surface $\partial_0E$, with $E$ on the concave side. It also has a boundary at infinity $\partial_\infty E$. Each connected component of the lift to $\tilde M=\HH^3$ of $E$ is a domain, say $\tilde E$, which is bounded in $\HH^3$ by a concave pleated surface $\tilde{\partial_0 E}$ which is a connected component of the lift of $\partial_0 E$, and has boundary at infinity by an open subset of $\partial_\infty\HH^3$ which is a connected component $\tilde{\partial_\infty E}$ of the lift of $\partial_\infty E$ to $\partial_\infty\HH^3$.

In this situation, the fundamental group of $E$, $\pi_1E$, acts properly discontinously on $\tilde E$, but also on $\tilde{\partial_\infty E}$. Considered in the projective model of $\HH^3$, this action extends to a properly discontinuous action to an open subset $\tilde{E^*}$ of $dS^3$, which can be defined as the set of points dual to the oriented totally geodesic planes contained in $\tilde E$. The quotient is a GHMC de Sitter spacetime $E^*$, which can be considered as a natural projective extension of $E$. By construction, $E$ and $E^*$ share the same boundary at infinity, which for $E^*$ is a boundary at future infinity.

This construction establishes a one-to-one correspondence between hyperbolic ends and GHMC de Sitter spacetimes, which also extends in higher dimension (see \cite{scannell,scannell:deformations}).

As for polyhedra, a smooth, closed, locally convex surface in $E$ has a dual in $E^*$, and the induced metric on  one is the third fundamental form on the other, a fact that is important in the questions considered in Section \ref{ssc:core} and in Section \ref{ssc:larger}. The dual volume of a region bounded by a surface can be interpreted as a volume bounded by the dual surface.

\subsection{Anti-de Sitter spacetimes as analogs of quasifuchsian manifolds}
\label{ssc:ads}

GHMC AdS spacetimes are quotients of a convex domain in $AdS^3$ by the fundamental group of a surface. Mess \cite{mess} discovered a remarkable analogy between GHMC AdS spacetimes and quasifuchsian hyperbolic manifolds. For instance:
\begin{itemize}
\item the moduli space $\cGH_S$ of GHMC AdS spacetimes on $S\times \R$, where $S$ is a closed surface of genus at least $2$, is parameterized by $\cT_S\times \cT_S$, as in for quasifuchsian hyperbolic manifolds through the Bers Simultaneous Uniformization Theorem. 
\item a GHMC AdS spacetime $M\in \cGH_S$ contains a smallest non-empty geodesically convex subset, its convex core $C(M)$. The boundary of $C(M)$ is (except for Fuchsian spacetimes) the disjoint union of two pleated surfaces homeomorphic to $S$, each equipped with a hyperbolic metric $m_\pm$ and pleated along a measured lamination $l_\pm$.
\item Mess conjectured --- following Thurston for quasifuchsian manifolds --- that any pair of hyperbolic metrics could be uniquely realized as $(m_-, m_+)$ for a certain $g\in \cGH_S$, and that any filling pair of measured laminations could be uniquely realized as $(l_-, l_+)$ for a certain $g\in \cGH_S$.
\end{itemize}

It was later proved that any pair or hyperbolic metric can be realized as the induced metric on $\partial C(M)$ \cite{diallo2013}, and that any filling pair of measured laminations can be realized as bending laminations \cite{earthquakes}. Moreover, any pair of smooth Riemannian metrics on $S$ of curvature $K<1$ can be realized as the induced metric on a ``larger'' geodesically convex subset of a GHMC AdS spacetime \cite{tamburelli2016}, which implies (through a duality between $AdS^3$ and itself) the same result for the third fundamental form. This leaves a number of questions open.

\begin{question}
  Given $M\in \cGH_S$:
  \begin{itemize}
  \item Is $M$ uniquely determined by the induced metric on $\partial C(M)$?
  \item Is $M$ uniquely determined by the measured bending lamination on $\partial C(M)$?
  \item If $N\subset M$ is a compact, geodesically convex subset with smooth boundary, is the pair $(M,N)$ uniquely determined by the induced metric on $\partial N$?
  \end{itemize}
\end{question}


\subsection{Minkowski and half-space geometries}
\label{ssc:minko}

GHMC Minkowski spacetimes are quotients of a convex domain in $\R^{2,1}$, which is future (or past) complete. Their holonomy representation is of the form $\rho: \pi_1S\to PSL(2,\R)\ltimes \R^{2,1}$. They occur in pairs, one future-convex and one past-convex, sharing the same holonomy. They can be obtained as first-order deformations of Fuchsian AdS or dS GHMC spacetimes, with the linear part of the holonomy, in $PSL(2,\R)$, corresponding to the Fuchsian holonomy, and the translation part, in $\R^{2,1}$, corresponding to the first-order deformation.

In this setting, there is an analog of Theorem \ref{tm:hmcb}: given two smooth metrics $h_-, h_+$ of curvature $K<0$ on a surface, there is a unique pair of GHMC Minkowski spacetimes (sharing the same holonomy) such that one is past-complete and contains a closed surface with induced metric $h_-$, while the other is future-complete and contains a closed surface with induced metric $h_+$, see \cite{smith2020weyl}.

GHMC Minkowski spacetimes are dual to the half-pipe manifolds introduced and studied by Danciger \cite{danciger:transition}. Those half-pipe manifolds are equipped with a degenerate metric, but they also have well-defined notion of convex core, with a boundary which has a hyperbolic induced metric and is pleated along a measured geodesic lamination (see e.g. \cite{fillastre-seppi}).


\bibliographystyle{amsplain}
\bibliography{/home/jean-marc/Dropbox/papiers/outils/biblio}

\def\cprime{$'$}
\providecommand{\bysame}{\leavevmode\hbox to3em{\hrulefill}\thinspace}
\providecommand{\MR}{\relax\ifhmode\unskip\space\fi MR }
\providecommand{\MRhref}[2]{%
  \href{http://www.ams.org/mathscinet-getitem?mr=#1}{#2}
}
\providecommand{\href}[2]{#2}
\begin{thebibliography}{10}

\bibitem{alex}
Alexander~D. Alexandrov, \emph{Convex polyhedra}, Springer Monographs in
  Mathematics, Springer-Verlag, Berlin, 2005, Translated from the 1950 Russian
  edition by N. S. Dairbekov, S. S. Kutateladze and A. B. Sossinsky, With
  comments and bibliography by V. A. Zalgaller and appendices by L. A. Shor and
  Yu. A. Volkov.

\bibitem{andreev}
E.M. Andreev, \emph{Convex polyhedra in {Lobacevskii} space}, Mat. Sb.(N.S.)
  \textbf{81 (123)} (1970), 445--478.

\bibitem{andreev-ideal}
\bysame, \emph{On convex polyhedra of finite volume in {Lobacevskii} space},
  Math. USSR Sbornik \textbf{12 (3)} (1971), 225--259.

\bibitem{bishop:explicit}
Christopher~J. Bishop, \emph{An explicit constant for {Sullivan}'s convex hull
  theorem}, In the tradition of Ahlfors and Bers, III. Proceedings of the 2nd
  Ahlfors-Bers colloquium, Storrs, CT, USA, October 18--21, 2001, Providence,
  RI: American Mathematical Society (AMS), 2004, pp.~41--69 (English).

\bibitem{bishop:weil-petersson}
\bysame, \emph{Weil-petersson curves, {{\(\beta\)}}-numbers, and minimal
  surfaces}, Ann. Math. (2) \textbf{202} (2025), no.~1, 111--188 (English).

\bibitem{bonahon}
Francis Bonahon, \emph{A {Schl\"afli}-type formula for convex cores of
  hyperbolic 3-manifolds}, J. Differential Geom. \textbf{50} (1998), no.~1,
  25--58.

\bibitem{bonahon-variations}
Francis Bonahon, \emph{Variations of the boundary geometry of {$3$}-dimensional
  hyperbolic convex cores}, J. Differential Geom. \textbf{50} (1998), no.~1,
  1--24. \MR{MR1678469 (2000j:57048)}

\bibitem{bonahon:laminations}
\bysame, \emph{Geodesic laminations on surfaces}, Laminations and foliations in
  dynamics, geometry and topology ({S}tony {B}rook, {NY}, 1998), Contemp.
  Math., vol. 269, Amer. Math. Soc., Providence, RI, 2001, pp.~1--37.
  \MR{1810534 (2001m:57023)}

\bibitem{bonahon-otal}
Francis Bonahon and Jean-Pierre Otal, \emph{Laminations mesur\'ees de plissage
  des vari\'et\'es hyperboliques de dimension 3}, Ann. Math. \textbf{160}
  (2004), 1013--1055.

\bibitem{convexhull}
Francesco Bonsante, Jeffrey Danciger, Sara Maloni, and Jean-Marc Schlenker,
  \emph{The induced metric on the boundary of the convex hull of a quasicircle
  in hyperbolic and anti--de {S}itter geometry}, Geom. Topol. \textbf{25}
  (2021), no.~6, 2827--2911, With an appendix by Boubacar Diallo. \MR{4347307}

\bibitem{earthquakes}
Francesco Bonsante and Jean-Marc Schlenker, \emph{Fixed points of compositions
  of earthquakes}, Duke Math. J. \textbf{161} (2012), no.~6, 1011--1054.
  \MR{2913100}

\bibitem{bridgeman-canary-yarmola}
M.~Bridgeman, R.~Canary, and A.~Yarmola, \emph{An improved bound for
  {Sullivan}'s convex hull theorem}, Proc. Lond. Math. Soc. (3) \textbf{112}
  (2016), no.~1, 146--168 (English).

\bibitem{bridgeman-brock-bromberg}
Martin Bridgeman, Jeffrey Brock, and Kenneth Bromberg, \emph{Schwarzian
  derivatives, projective structures, and the {W}eil-{P}etersson gradient flow
  for renormalized volume}, Duke Math. J. \textbf{168} (2019), no.~5, 867--896.
  \MR{3934591}

\bibitem{bridgeman-brock-bromberg:gradient}
\bysame, \emph{The {W}eil-{P}etersson gradient flow of renormalized volume and
  3-dimensional convex cores}, Geom. Topol. \textbf{27} (2023), no.~8,
  3183--3228. \MR{4668096}

\bibitem{bridgeman-bromberg-pallete}
Martin Bridgeman, Kenneth Bromberg, and Franco Vargas~Pallete, \emph{The
  {Weil--Petersson} gradient flow of renormalized volume on a {Bers} slice has
  a global attracting fixed point}, arXiv preprint arXiv:2105.01207 (2021).

\bibitem{bridgeman-bromberg-pallete:convergence}
\bysame, \emph{Convergence of the gradient flow of renormalized volume to
  convex cores with totally geodesic boundary}, Compos. Math. \textbf{159}
  (2023), no.~4, 830--859. \MR{4571574}

\bibitem{BBVPW}
Martin Bridgeman, Kenneth Bromberg, Franco {Vargas Pallete}, and Yilin Wang,
  \emph{Universal {Liouville} action as a renormalized volume and its gradient
  flow}, 2025.

\bibitem{bridgeman-canary:renormalized}
Martin Bridgeman and Richard~D. Canary, \emph{Renormalized volume and the
  volume of the convex core}, Ann. Inst. Fourier (Grenoble) \textbf{67} (2017),
  no.~5, 2083--2098. \MR{3732685}

\bibitem{brock:2003}
Jeffrey~F. Brock, \emph{The {W}eil-{P}etersson metric and volumes of
  3-dimensional hyperbolic convex cores}, J. Amer. Math. Soc. \textbf{16}
  (2003), no.~3, 495--535 (electronic). \MR{1969203 (2004c:32027)}

\bibitem{brock-bromberg:inflexibility2}
Jeffrey~F. Brock and Kenneth~W. Bromberg, \emph{Inflexibility,
  {W}eil-{P}etersson distance, and volumes of fibered 3-manifolds}, Math. Res.
  Lett. \textbf{23} (2016), no.~3, 649--674. \MR{3533189}

\bibitem{canary:conformal}
R.~D. Canary, \emph{The conformal boundary and the boundary of the convex
  core}, Duke Math. J. \textbf{106} (2001), no.~1, 193--207. \MR{1810370}

\bibitem{chen2022geometric}
Qiyu Chen and Jean-Marc Schlenker, \emph{The geometric data on the boundary of
  convex subsets of hyperbolic manifolds}, 2022, To appear, {\em Trans. Amer.
  Math. Soc.}

\bibitem{choudhury-markovic}
Diptaishik Choudhury and Vladimir Markovic, \emph{Measured foliations at
  infinity of quasi-fuchsian manifolds}, 2025.

\bibitem{filling}
Tommaso Cremaschi, Viola Giovannini, and Jean-Marc Schlenker, \emph{Filling
  {Riemann} surfaces by hyperbolic schottky manifolds of negative volume},
  2024, To appear, {\it J. Geom. Physics}.

\bibitem{danciger:transition}
Jeffrey Danciger, \emph{A geometric transition from hyperbolic to anti-de
  {S}itter geometry}, Geom. Topol. \textbf{17} (2013), no.~5, 3077--3134.
  \MR{3190306}

\bibitem{diallo2013}
Boubacar Diallo, \emph{Prescribing metrics on the boundary of convex cores of
  globally hyperbolic maximal compact ads 3-manifolds}, arXiv preprint
  arXiv:1303.7406 (2013), Published as the appendix in {\it The induced metric
  on the boundary of the convex hull of a quasicircle in hyperbolic and
  anti–de Sitter geometry, F. Bonsante, J. Danciger, S. Maloni, and J.-M.
  Schlenker, Geom. Topol., 25 (2021), pp. 2827--2911}.

\bibitem{uniqueness}
Bruno Dular and Jean-Marc Schlenker, \emph{Convex co-compact hyperbolic
  manifolds are determined by their pleating lamination}, 2024.

\bibitem{dumas-survey}
Emily Dumas, \emph{Complex projective structures}, Handbook of {T}eichm\"uller
  theory. {V}ol. {II}, IRMA Lect. Math. Theor. Phys., vol.~13, Eur. Math. Soc.,
  Z\"urich, 2008, pp.~455--508.

\bibitem{epstein:envelopes}
Charles~L Epstein, \emph{Envelopes of horospheres and {Weingarten} surfaces in
  hyperbolic 3-space}, arxiv:2401.12115 (2024), Original preprint from 1984.

\bibitem{epstein-marden}
D.~B.~A. Epstein and A.~Marden, \emph{Convex hulls in hyperbolic spaces, a
  theorem of {Sullivan}, and measured pleated surfaces}, Analytical and
  geometric aspects of hyperbolic space (D.~B.~A. Epstein, ed.), L.M.S. Lecture
  Note Series, vol. 111, Cambridge University Press, 1986.

\bibitem{epstein-marden-markovic}
D.~B.~A. Epstein, A.~Marden, and V.~Markovic, \emph{Quasiconformal
  homeomorphisms and the convex hull boundary}, Ann. of Math. (2) \textbf{159}
  (2004), no.~1, 305--336. \MR{2052356}

\bibitem{FLP}
A.~Fathi, F.~Laudenbach, and V.~Poenaru, \emph{Travaux de {T}hurston sur les
  surfaces}, Soci\'et\'e Math\'ematique de France, Paris, 1991, S\'eminaire
  Orsay, Reprint of {\it Travaux de Thurston sur les surfaces}, Soc.\ Math.\
  France, Paris, 1979 [MR 82m:57003], Ast\'erisque No. 66-67 (1991).

\bibitem{fillastre-seppi}
Fran{\c{c}}ois Fillastre and Andrea Seppi, \emph{Spherical, hyperbolic and
  other projective geometries: convexity, duality, transitions}, arXiv preprint
  arXiv:1611.01065 (2016).

\bibitem{giovannini:adapted}
Viola Giovannini, \emph{Adapted renormalized volume for hyperbolic 3-manifolds
  with compressible boundary}, 2025.

\bibitem{graham}
C.~Robin Graham, \emph{Volume and area renormalizations for conformally compact
  {E}instein metrics}, The {P}roceedings of the 19th {W}inter {S}chool
  ``{G}eometry and {P}hysics'' ({S}rn\'\i, 1999), no.~63, 2000, pp.~31--42.
  \MR{MR1758076 (2002c:53073)}

\bibitem{graham-witten}
C.~Robin Graham and Edward Witten, \emph{Conformal anomaly of submanifold
  observables in {A}d{S}/{CFT} correspondence}, Nuclear Phys. B \textbf{546}
  (1999), no.~1-2, 52--64. \MR{1682674 (2000h:81286)}

\bibitem{guillarmou-moroianu-rochon}
Colin Guillarmou, Sergiu Moroianu, and Fr\'ed\'eric Rochon, \emph{Renormalized
  volume on the {T}eichm\"uller space of punctured surfaces}, Ann. Sc. Norm.
  Super. Pisa Cl. Sci. (5) \textbf{17} (2017), no.~1, 323--384. \MR{3676051}

\bibitem{Skenderis}
Mark Henningson and Kostas Skenderis, \emph{The holographic {Weyl} anomaly},
  JHEP \textbf{9807:023} (1998).

\bibitem{HR}
Craig~D. Hodgson and Igor Rivin, \emph{A characterization of compact convex
  polyhedra in hyperbolic 3-space}, Invent. Math. \textbf{111} (1993), 77--111.

\bibitem{kerckhoff:analytic}
Steven~P. Kerckhoff, \emph{Earthquakes are analytic}, Comment. Math. Helv.
  \textbf{60} (1985), no.~1, 17--30. \MR{787659 (86m:57014)}

\bibitem{kojima-mcshane}
Sadayoshi Kojima and Greg McShane, \emph{Normalized entropy versus volume for
  pseudo-{A}nosovs}, Geom. Topol. \textbf{22} (2018), no.~4, 2403--2426.
  \MR{3784525}

\bibitem{Holography}
Kirill Krasnov, \emph{Holography and {R}iemann surfaces}, Adv. Theor. Math.
  Phys. \textbf{4} (2000), no.~4, 929--979. \MR{MR1867510 (2002k:81230)}

\bibitem{volume}
Kirill Krasnov and Jean-Marc Schlenker, \emph{On the renormalized volume of
  hyperbolic 3-manifolds}, Comm. Math. Phys. \textbf{279} (2008), no.~3,
  637--668. \MR{MR2386723}

\bibitem{L4}
Fran\c{c}ois Labourie, \emph{M\'etriques prescrites sur le bord des
  vari\'et\'es hyperboliques de dimension 3}, J. Differential Geom. \textbf{35}
  (1992), 609--626.

\bibitem{L6}
Fran{\c{c}}ois Labourie, \emph{Probl\`eme de {M}inkowski et surfaces \`a
  courbure constante dans les vari\'et\'es hyperboliques}, Bull. Soc. Math.
  France \textbf{119} (1991), no.~3, 307--325. \MR{92k:53077}

\bibitem{lecuire}
Cyril Lecuire, \emph{Plissage des vari\'et\'es hyperboliques de dimension 3},
  Invent. Math. \textbf{164} (2006), no.~1, 85--141. \MR{MR2207784
  (2006m:57029)}

\bibitem{maldacena}
J.~Maldacena, \emph{The large {N} limit of superconformal field theories and
  supergravity}, Adv. Theor. Math. Phys. \textbf{2} (1998), 231--252.

\bibitem{marden:hyperbolic}
Albert Marden, \emph{Hyperbolic manifolds}, Cambridge University Press,
  Cambridge, 2016, An introduction in 2 and 3 dimensions. \MR{3586015}

\bibitem{matsuzaki-taniguchi}
Katsuhiko Matsuzaki and Masahiko Taniguchi, \emph{Hyperbolic manifolds and
  {K}leinian groups}, Oxford Mathematical Monographs, The Clarendon Press,
  Oxford University Press, New York, 1998, Oxford Science Publications.
  \MR{1638795}

\bibitem{mazzoli2019-1}
Filippo Mazzoli, \emph{The dual volume of quasi-{F}uchsian manifolds and the
  {W}eil-{P}etersson distance}, Trans. Amer. Math. Soc. \textbf{375} (2022),
  no.~1, 695--723. \MR{4358680}

\bibitem{McMullen}
Curtis~T. McMullen, \emph{The moduli space of {R}iemann surfaces is {K}\"ahler
  hyperbolic}, Ann. of Math. (2) \textbf{151} (2000), no.~1, 327--357.
  \MR{MR1745010 (2001m:32032)}

\bibitem{mesbah2024}
Abderrahim Mesbah, \emph{The prescribed metric on the boundary of convex
  subsets of anti-de sitter space with a quasi-circle as ideal boundary}, 2024,
  To appear, {\it J. Geom. Anal.}

\bibitem{mess}
Geoffrey Mess, \emph{Lorentz spacetimes of constant curvature}, Geom. Dedicata
  \textbf{126} (2007), 3--45. \MR{MR2328921}

\bibitem{milnor-schlafli}
John Milnor, \emph{The {Schl\"afli} differential equality}, Collected papers,
  vol. 1, Publish or Perish, 1994.

\bibitem{nehari-bams}
Zeev Nehari, \emph{The {S}chwarzian derivative and schlicht functions}, Bull.
  Amer. Math. Soc. \textbf{55} (1949), 545--551. \MR{0029999}

\bibitem{otal-hyperbolisation}
Jean-Pierre Otal, \emph{Le th\'eor\`eme d'hyperbolisation pour les vari\'et\'es
  fibr\'ees de dimension 3}, Ast\'erisque \textbf{235} (1996), x+159.

\bibitem{VPWW}
Franco~Vargas Pallete, Yilin Wang, and Catherine Wolfram, \emph{Epstein curves
  and holography of the {Schwarzian} action}, 2025.

\bibitem{Po}
Aleksei~V. Pogorelov, \emph{Extrinsic geometry of convex surfaces}, American
  Mathematical Society, 1973, Translations of Mathematical Monographs. Vol. 35.

\bibitem{prosanov:dual}
Roman Prosanov, \emph{Dual metrics on the boundary of strictly polyhedral
  hyperbolic 3-manifolds}, 2022, arxiv:2203.16971. To appear, {\em J.
  Differential Geom.}

\bibitem{prosanov:polyhedral}
\bysame, \emph{Hyperbolic 3-manifolds with boundary of polyhedral type}, 2022,
  arxiv:2210.17271.

\bibitem{rivin-annals}
Igor Rivin, \emph{A characterization of ideal polyhedra in hyperbolic 3-space},
  Annals of Math. \textbf{143} (1996), 51--70.

\bibitem{sem-era}
Igor Rivin and Jean-Marc Schlenker, \emph{The {Schl\"afli} formula in
  {Einstein} manifolds with boundary}, Electronic Research Announcements of the
  A.M.S. \textbf{5} (1999), 18--23.

\bibitem{rohde-wang}
Steffen Rohde and Yilin Wang, \emph{The {Loewner} energy of loops and
  regularity of driving functions}, Int. Math. Res. Not. \textbf{2021} (2021),
  no.~10, 7433--7469 (English).

\bibitem{scannell}
Kevin~P. Scannell, \emph{Flat conformal structures and the classification of de
  {S}itter manifolds}, Comm. Anal. Geom. \textbf{7} (1999), no.~2, 325--345.
  \MR{1685590 (2000a:53122)}

\bibitem{scannell:deformations}
\bysame, \emph{Infinitesimal deformations of some {${\rm SO}(3,1)$} lattices},
  Pacific J. Math. \textbf{194} (2000), no.~2, 455--464. \MR{1760793
  (2001c:57018)}

\bibitem{schlafli1858}
Ludwig Schl\"afli, \emph{Note sur quelques propri\'et\'es des poly\`edres},
  Journal de Math\'ematiques Pures et Appliqu\'ees \textbf{3} (1858), 269--301.

\bibitem{hmcb}
Jean-Marc Schlenker, \emph{Hyperbolic manifolds with convex boundary}, Invent.
  Math. \textbf{163} (2006), no.~1, 109--169. \MR{MR2208419 (2006m:57023)}

\bibitem{compare}
\bysame, \emph{The renormalized volume and the volume of the convex core of
  quasifuchsian manifolds}, Math. Res. Lett. \textbf{20} (2013), no.~4,
  773--786, Corrected version available as arXiv:1109.6663v4. \MR{3188032}

\bibitem{volumes}
\bysame, \emph{Volumes of quasifuchsian manifolds}, Surveys in differential
  geometry 2020. {S}urveys in 3-manifold topology and geometry, Surv. Differ.
  Geom., vol.~25, Int. Press, Boston, MA, 2022, pp.~319--353. \MR{4479756}

\bibitem{averages}
Jean-Marc Schlenker and Edward Witten, \emph{No ensemble averaging below the
  black hole threshold}, J. High Energy Phys. (2022), no.~7, Paper No. 143, 50.
  \MR{4458319}

\bibitem{shen:weil-petersson}
Yuliang Shen, \emph{Weil--{P}etersson {T}eichm\"{u}ller space}, Amer. J. Math.
  \textbf{140} (2018), no.~4, 1041--1074. \MR{3828040}

\bibitem{slutskiy:compact}
Dmitriy Slutskiy, \emph{Compact domains with prescribed convex boundary metrics
  in quasi-fuchsian manifolds}, Bull. Soc. Math. France \textbf{146} (2018),
  no.~2, 309--353.

\bibitem{smith2020weyl}
Graham Smith, \emph{On the {Weyl} problem in {Minkowski} space}, arXiv preprint
  arXiv:2005.01137 (2020).

\bibitem{sullivan:travaux}
Dennis Sullivan, \emph{Travaux de {Thurston} sur les groupes quasi-{Fuchsiens}
  et les vari{\'e}t{\'e}s hyperboliques de dimension 3 fibrees sur {{\({S}^
  1\)}}.}, S{\'e}min. {Bourbaki}, 32e ann{\'e}e, {Vol}. 1979/80, {Exp}. 554,
  {Lect}. {Notes} {Math}. 842, 196-214 (1981)., 1981.

\bibitem{TZ-schottky}
L.~Takhtajan and P.~Zograf, \emph{On uniformization of {Riemann} surfaces and
  the {Weil--Petersson} metric on the {Teichm\"uller} and {Schottky} spaces},
  Mat. Sb. \textbf{132} (1987), 303--320, English translation in {\it Math.
  USSR Sb. 60:297-313, 1988}.

\bibitem{takhtajan-teo}
Leon~A. Takhtajan and Lee-Peng Teo, \emph{Liouville action and
  {W}eil-{P}etersson metric on deformation spaces, global {K}leinian
  reciprocity and holography}, Comm. Math. Phys. \textbf{239} (2003), no.~1-2,
  183--240. \MR{MR1997440 (2005c:32021)}

\bibitem{takhtajan-teo:universal}
\bysame, \emph{Weil--{P}etersson metric on the universal {T}eichm\"uller
  space}, Mem. Amer. Math. Soc. \textbf{183} (2006), no.~861, viii+119.
  \MR{2251887 (2007e:32011)}

\bibitem{tamburelli2016}
Andrea Tamburelli, \emph{Prescribing metrics on the boundary of anti--de
  {S}itter 3-manifolds}, Int. Math. Res. Not. IMRN (2018), no.~5, 1281--1313.
  \MR{3801462}

\bibitem{thurston-notes}
William~P. Thurston, \emph{Three-dimensional geometry and topology.},
  Originally notes of lectures at Princeton University, 1979. Recent version
  available on http://www.msri.org/publications/books/gt3m/, 1980.

\bibitem{pallete:schottky}
Franco Vargas~Pallete, \emph{Upper bounds on renormalized volume for {Schottky}
  groups}, Math. Ann. \textbf{392} (2025), no.~1, 733--750 (English).

\bibitem{wang:equivalent}
Yilin Wang, \emph{Equivalent descriptions of the {Loewner} energy}, Invent.
  Math. \textbf{218} (2019), no.~2, 573--621 (English).

\bibitem{witten}
E.~Witten, \emph{Anti de {Sitter} space and holography}, Adv. Theor. Math.
  Phys. \textbf{2} (1998), 253--291.

\end{thebibliography}
\end{document}